\documentclass[A4,twoside,11pt,leqno]{article}
\usepackage{amsfonts}
\usepackage{amssymb}
\usepackage{amsthm}
\usepackage[tbtags]{amsmath}
\usepackage{doc}
\usepackage{latexsym}
\usepackage{amscd}
\usepackage[frame,cmtip,arrow,matrix,line,graph,curve]{xy}
\usepackage{graphics}
\usepackage{epsfig}
\usepackage{eucal}
\usepackage{mathrsfs}

\font\headd=cmr8
\begin{document}

\thispagestyle{plain} \markboth{}{} {\small {\addtocounter{page}{0} %
\noindent {\large \textbf{New Generalized Fixed Point Results on $S_{b}$%
-Metric Spaces}} \footnote{%
{}\newline
\\[-0.7cm]
* Corresponding Author.\newline
Received March 00, 2014; revised May 00, 2014; accepted November 00, 2014.%
} \vspace{0.15in}\newline
\noindent \textsc{N\.{I}HAL TA\c{S}$^{\ast }$ and N\.{I}HAL YILMAZ \"{O}ZG%
\"{U}R} \newline
\textit{Department of Mathematics, Bal\i kesir University, Bal\i kesir
10145, Turkey\newline
e-mail} : {\verb|nihaltas@balikesir.edu.tr|} \textit{and} {%
\verb|nihal@balikesir.edu.tr|} \vspace{0.15in}\newline
{\footnotesize \textsc{Abstract. }Recently $S_{b}$-metric spaces have been
introduced as the generalizations of metric and }}}${\small {{\footnotesize {%
S}}}}${\small {{\footnotesize {-metric spaces. In this paper} we investigate
some basic properties of this new space. We generalize the classical Banach's contraction
principle using the theory of a complete $S_{b}$-metric space. Also we give an application to linear equation systems using
the $S_{b}$-metric which is generated by a metric. } \vspace{0.2in}\newline
\pagestyle{myheadings}
\markboth{\headd N. TA\c{S} and N. YILMAZ \"{O}ZG\"{U}R$~~~~~~~~~~~~~~~~~~~~~~~~~~~~~~~~~~~~~~~~~~~~~\,$}
 {\headd $~~~~~~~~~~~~~~~~~$New Generalized Fixed Point Results on $S_{b}$%
-Metric Spaces}%
\newline
\noindent \textbf{1. Introduction and mathematical preliminaries} %
\setcounter{equation}{0} \renewcommand{\theequation}{1.\arabic{equation}}
\vspace{0.1in}\newline
\indent\label{intro} Metric spaces and fixed point theorems are very
important in many areas of mathematics. Some generalizations of metric
spaces and fixed points of various contractive mappings have been studied
extensively. }Bakhtin introduced $b$-metric spaces as a generalization of
metric spaces \cite{Bakhtin-1989}. Mustafa and Sims defined the concept of a
generalized metric space which is called a $G$-metric space \cite%
{Mustafa-2006}. Sedghi, Shobe and Aliouche gave the notion of an $S$-metric
space and proved some fixed-point theorems for a self-mapping on a complete $%
S$-metric space \cite{Sedghi-2012}. Aghajani, Abbas and Roshan presented a
new type of metric which is called $G_{b}$-metric and studied some
properties of this metric \cite{Aghajani-2014-1}. Since then, many authors
obtained several fixed-point results in the various generalized metric
spaces (see \cite{Aghajani-2014-2}, \cite{An-2013}, \cite{Dung-2014}, \cite%
{Ege-2016}, \cite{Hieu-2015}, \cite{Hussain-2015}, \cite{Mohanta-2012}, \cite%
{Nihal-2015-1}, \cite{Nihal-2015-2}, \cite{Nihal-2015-3}, \cite%
{Sedghi-2014-1} and \cite{Sedghi-2014-2} for more details). Also some
applications of fixed point theory were studied on various metric spaces. It
was given several applications of the Banach's contraction principle in many
areas such as integral equations, linear equations, differential equations
etc. For example, the present authors were investigated complex and
differential applications on $S$-metric spaces (see \cite{Nihal-complex} and
\cite{Nihal-picard}). }

{\small Recently, the concept of an $S_{b}$-metric space as a generalization
of metric spaces and $S$-metric spaces have been introduced in \cite%
{Sedghi-2016} and a common fixed point theorem for four mappings have been
studied on a complete $S_{b}$-metric space. Also the notion of an $S_{b}$%
-metric was generalized to the notion of an $A_{b}$-metric in \cite%
{Ughade-2016}. An $S_{b}$-metric space is a space with three dimensions
whereas an $A_{b}$-metric space is a space with $n$ dimensions. When $n=3$,
these two notions of \textquotedblleft an $S_{b}$-metric\textquotedblright\
and \textquotedblleft an }${\small A}${\small $_{b}$-metric%
\textquotedblright\ coincide. Some fixed point theorems were given under
different contraction and expansion type conditions (see \cite{Ughade-2016}
for more details). }

{\small In this paper we consider a complete $S_{b}$-metric space and prove
two generalizations of the classical Banach's fixed point result. In Section
$2$, we recall some known definitions. In Section $3$, we deal with the
notion of an $S_{b}$-metric and investigate some properties of $S_{b}$%
-metric spaces. We study some relationships between an $S_{b}$-metric and
other some metrics. In Section $4$, we prove the Banach's contraction
principle on a complete $S_{b}$-metric space and give a new fixed point
theorem as a generalization of the Banach's contraction principle with a
counterexample. In Section $5$, we present an application to linear
equations on an $S_{b}$-metric space $(X,S_{1})$. }

{\small \setcounter{equation}{0} \renewcommand{\theequation}{2.%
\arabic{equation}} \noindent \textbf{2. Preliminaries}\vspace{0.1in}\newline
\indent\label{sec:1} In this section we recall the following definitions. }

{\small \textbf{Definition 1.1.} \label{def1} \cite{Bakhtin-1989} Let $X$ be
a nonempty set, $b\geq 1$ be a given real number and $d:X\times X\rightarrow
\lbrack 0,\infty )$ be a function satisfying the following conditions for
all $x,y,z\in X$. }

{\small $(b1)$ $d(x,y)=0$ if and only if $x=y$. }

{\small $(b2)$ $d(x,y)=d(y,x)$. }

{\small $(b3)$ $d(x,z)\leq b[d(x,y)+d(y,z)]$. }

{\small Then the function $d$ is called a $b$-metric on $X$ and the pair $%
(X,d)$ is called a $b$-metric space. }

{\small \textbf{Definition 1.2.} \label{def2} \cite{Mustafa-2006} Let $X$ be
a nonempty set and $G:X\times X\times X\rightarrow \lbrack 0,\infty )$ be a
function satisfying the following conditions. }

{\small $(G1)$ $G(x,y,z)=0$ if $x=y=z$. }

{\small $(G2)$ $0<G(x,x,y)$ for all $x,y\in X$ with $x\neq y$. }

{\small $(G3)$ $G(x,x,y)\leq G(x,y,z)$ for all $x,y,z\in X$ with $y\neq z$. }

{\small $(G4)$ $G(x,y,z)=G(x,z,y)=G(y,z,x)=\cdots $. }

{\small $(G5)$ $G(x,y,z)\leq G(x,a,a)+G(a,y,z)$ for all $x,y,z,a\in X$. }

{\small Then the function $G$ is called a generalized metric or a $G$-metric
on $X$ and the pair $(X,G)$ is called a $G$-metric space. }

{\small \textbf{Definition 1.3.} \label{def3} \cite{Aghajani-2014-1} Let $X$
be a nonempty set, $b\geq 1$ be a given real number and $G_{b}:X\times
X\times X\rightarrow \lbrack 0,\infty )$ be a function satisfying the
following conditions. }

{\small $(G_{b}1)$ $G_{b}(x,y,z)=0$ if $x=y=z$. }

{\small $(G_{b}2)$ $0<G_{b}(x,x,y)$ for all $x,y\in X$ with $x\neq y$. }

{\small $(G_{b}3)$ $G_{b}(x,x,y)\leq G_{b}(x,y,z)$ for all $x,y,z\in X$ with
$y\neq z $. }

{\small $(G_{b}4)$ $G_{b}(x,y,z)=G_{b}(x,z,y)=G_{b}(y,z,x)=\cdots $. }

{\small $(G_{b}5)$ $G_{b}(x,y,z)\leq b[G_{b}(x,a,a)+G_{b}(a,y,z)]$ for all $%
x,y,z,a\in X$. }

{\small Then the function $G_{b}$ is called a generalized $b$-metric or a $%
G_{b}$-metric on $X$ and the pair $(X,G_{b})$ is called a $G_{b}$-metric
space. }

{\small \textbf{Definition 1.4.} \label{def4} \cite{Sedghi-2012} Let $X$ be
a nonempty set and $S:X\times X\times X\rightarrow \lbrack 0,\infty )$ be a
function satisfying the following conditions for all $x,y,z,a\in X$. }

{\small $(S1)$ $S(x,y,z)=0$ if and only if $x=y=z$. }

{\small $(S2)$ $S(x,y,z)\leq S(x,x,a)+S(y,y,a)+S(z,z,a)$. }

{\small Then the function $S$ is called an $S$-metric on $X$ and the pair $%
(X,S)$ is called an $S$-metric space. }

{\small We use the following lemma in the next sections. }

{\small \textbf{Lemma 1.1.} \label{lem6} \cite{Sedghi-2012} Let $(X,S)$ be
an $S$-metric space. Then we have%
\begin{equation*}
S(x,x,y)=S(y,y,x)\text{.}
\end{equation*}%
\vspace{0.2in}\newline
\setcounter{equation}{0} \renewcommand{\theequation}{2.\arabic{equation}}
\noindent \textbf{3. $S_{b}$-Metric spaces} \vspace{0.1in}\newline
\indent\label{sec:2} In this section we recall the notion of an $S_{b}$%
-metric space and study some properties of this space. }

{\small \textbf{Definition 2.1.} \cite{Sedghi-2016}\label{def5} Let $X$ be a
nonempty set and $b\geq 1$ be a given real number. A function $S_{b}:X\times
X\times X\rightarrow \lbrack 0,\infty )$ is said to be $S_{b}$-metric if and
only if for all $x,y,z,a\in X$ the following conditions are satisfied: }

{\small $(S_{b}1)$ $S_{b}(x,y,z)=0$ if and only if $x=y=z$, }

{\small $(S_{b}2)$ $S_{b}(x,y,z)\leq
b[S_{b}(x,x,a)+S_{b}(y,y,a)+S_{b}(z,z,a)]$. }

{\small The pair $(X,S_{b})$ is called an $S_{b}$-metric space. }

{\small We note that $S_{b}$-metric spaces are the generalizations of $S$%
-metric spaces since every $S$-metric is an $S_{b}$-metric with $b=1$. But
the converse statement is not always true (see \cite{Sedghi-2016} for more
details). In the following we give another example of an $S_{b}$-metric
which is not an $S$-metric on $X$. }

{\small \textbf{Example 2.1.} \label{exm1} Let $X=%
\mathbb{R}
$ and the function $S_{b}$ be defined as%
\begin{equation*}
S_{b}(x,y,z)=\dfrac{1}{16}(\left\vert x-y\right\vert +\left\vert
y-z\right\vert +\left\vert x-z\right\vert )^{2}\text{.}
\end{equation*}%
Then the function $S_{b}$ is an $S_{b}$-metric with $b=4$, but it is not an $%
S$-metric. Indeed, for $x=4$, $y=6$, $z=8$ and $a=5$ we get%
\begin{equation*}
S_{b}(4,6,8)=4\text{, }S_{b}(4,4,5)=\dfrac{1}{4}\text{, }S_{b}(6,6,5)=\dfrac{%
1}{4}\text{, }S_{b}(8,8,5)=\dfrac{9}{4}\text{.}
\end{equation*}%
Hence we have%
\begin{equation*}
S_{b}(4,6,8)=4\leq S_{b}(4,4,5)+S_{b}(6,6,5)+S_{b}(8,8,5)=\dfrac{11}{4}\text{%
,}
\end{equation*}%
which is a contradiction with $(S2)$. }

{\small \textbf{Definition 2.2.} \label{def9} Let $(X,S_{b})$ be an $S_{b}$%
-metric space and $b>1$. An $S_{b}$-metric $S_{b}$ is called symmetric if
\begin{equation}
S_{b}(x,x,y)=S_{b}(y,y,x)\text{,}  \label{symmetry}
\end{equation}%
for all $x,y\in X$. }

{\small In \cite{Mlaiki-2016}, it was given a definition of an $S_{b}$%
-metric with the symmetry condition \textquotedblleft $%
S_{b}(x,x,y)=S_{b}(y,y,x)$\textquotedblright\ (see Definition 1.3 on page
132). However in the definition of an $S_{b}$-metric the symmetry condition (%
\ref{symmetry}) is not necessary. In fact, for $b=1$ the $S_{b}$-metric
induced to an $S$-metric. It is known that the symmetry condition (\ref%
{symmetry}) is automatically satisfied by an $S$-metric (see Lemma 2.5 on
page 260 in \cite{Sedghi-2012}). So Definition 2.1 of an $S_{b}$-metric is
more general than given in \cite{Mlaiki-2016}. }

{\small We give the following examples of symmetric $S_{b}$-metric and
non-symmetric $S_{b}$-metric, respectively. }

{\small \textbf{Example 2.2.} \label{exm7} Let $(X,d)$ be a metric space and
the function $S_{b}:X\times X\times X\rightarrow \lbrack 0,\infty )$ be
defined as%
\begin{equation*}
S_{b}(x,y,z)=\left[ d(x,y)+d(y,z)+d(x,z)\right] ^{p}\text{,}
\end{equation*}%
for all $x,y,z\in X$ and $p>1$. Then it can be easily seen that $S_{b}$ is
an $S_{b}$-metric on $X$. Also the function $S_{b}$ satisfies the symmetric
condition $($\ref{symmetry}$)$. }

{\small \textbf{Example 2.3.} \label{exm8} Let $X=%
\mathbb{R}
$ and the function $S_{b}:X\times X\times X\rightarrow \lbrack 0,\infty )$
be defined as%
\begin{equation*}
\begin{array}{l}
S_{b}(0,0,1)=2\text{,} \\
S_{b}(1,1,0)=4\text{,} \\
S_{b}(x,y,z)=0\text{ if }x=y=z\text{,} \\
S_{b}(x,y,z)=1\text{ otherwise,}%
\end{array}%
\end{equation*}%
for all $x,y,z\in
\mathbb{R}
$. Then the function $S_{b}$ is an $S_{b}$-metric with $b\geq 2$ which is
not symmetric. }

{\small We define some topological concepts in the following: }

{\small \textbf{Definition 2.3.} \label{def10} Let $(X,S_{b})$ be an $S_{b}$%
-metric space, $x\in X$ and $A,B\subset X$. }

\begin{enumerate}
\item {\small We define the distance between the sets $A$ and $B$ by%
\begin{equation*}
S_{b}(A,A,B)=\inf \{S_{b}(x,x,y):x\in A,y\in B\}\text{.}
\end{equation*}
}

\item {\small We define the distance of the the point $x$ to the set $A$ by%
\begin{equation*}
S_{b}(x,x,A)=\inf \{S_{b}(x,x,y):y\in A\}\text{.}
\end{equation*}
}

\item {\small We define the diameter of $A$ by%
\begin{equation*}
\delta (A)=\sup \{S_{b}(x,x,y):x,y\in A\}\text{.}
\end{equation*}
}
\end{enumerate}

{\small Now we recall the definition of an open ball and a closed ball on $%
S_{b}$-metric spaces, respectively. }

{\small \textbf{Definition 2.4.}\label{def6} \cite{Sedghi-2016} Let $%
(X,S_{b}) $ be an $S_{b}$-metric space. The open ball $B_{S}^{b}(x,r)$ and
the closed ball $B_{S}^{b}[x,r]$ with a center $x$ and a radius $r$ are
defined by%
\begin{equation*}
B_{S}^{b}(x,r)=\{y\in X:S_{b}(y,y,x)<r\}
\end{equation*}%
and%
\begin{equation*}
B_{S}^{b}[x,r]=\{y\in X:S_{b}(y,y,x)\leq r\}\text{,}
\end{equation*}%
for $r>0$, $x\in X$, respectively. }

{\small \textbf{Example 2.4.} \label{exm2} Let us consider the $S_{b}$%
-metric space defined in Example $2.1$ as follows:%
\begin{equation*}
S_{b}(x,y,z)=\dfrac{1}{16}(\left\vert x-y\right\vert +\left\vert
y-z\right\vert +\left\vert x-z\right\vert )^{2}\text{,}
\end{equation*}%
for all $x,y,z\in
\mathbb{R}
$. Then we get%
\begin{equation*}
B_{S}^{b}(0,2)=\{y\in
\mathbb{R}
:S_{b}(y,y,0)<2\}=(-2\sqrt{2},2\sqrt{2})
\end{equation*}%
and%
\begin{equation*}
B_{S}^{b}[0,2]=\{y\in
\mathbb{R}
:S_{b}(y,y,0)\leq 2\}=[-2\sqrt{2},2\sqrt{2}]\text{.}
\end{equation*}
}

{\small \textbf{Definition 2.5.} \label{def7} Let $(X,S_{b})$ be an $S_{b}$%
-metric space and $X^{\prime }\subset X$. }

\begin{enumerate}
\item {\small If there exists $r>0$ such that $B_{S}^{b}(x,r)\subset
X^{\prime }$ for every $x\in X^{\prime }$ then $X^{\prime }$ is called an
open subset of $X$. }

\item {\small Let $\tau $ be the set of all $X^{\prime }\subset X$ with $%
x\in X^{\prime }$ such that there exists $r>0$ satisfying $%
B_{S}^{b}(x,r)\subset X^{\prime }$. Then $\tau $ is called the topology
induced by the $S_{b}$-metric. }

\item {\small $X^{\prime }$ is called $S_{b}$-bounded if there exists $r>0$
such that $S_{b}(x,x,y)<r$ for all $x,y\in X^{\prime }$. If $X^{\prime }$ is
$S_{b}$-bounded then we will write $\delta (X^{\prime })<\infty $. }
\end{enumerate}

{\small \textbf{Definition 2.6.} \label{def8}\cite{Sedghi-2016} Let $%
(X,S_{b})$ be an $S_{b}$-metric space. }

\begin{enumerate}
\item {\small A sequence $\{x_{n}\}$ in $X$ converges to $x$ if and only if $%
S_{b}(x_{n},x_{n},x)\rightarrow 0$ as $n\rightarrow \infty $, that is, for
each $\varepsilon >0$ there exists $n_{0}\in
\mathbb{N}
$ such that for all $n\geq n_{0}$, $S_{b}(x_{n},x_{n},x)<\varepsilon $. It
is denoted by%
\begin{equation*}
\underset{n\rightarrow \infty }{\lim }x_{n}=x\text{.}
\end{equation*}
}

\item {\small A sequence $\{x_{n}\}$ in $X$ is called a Cauchy sequence if
for each $\varepsilon >0$ there exists $n_{0}\in
\mathbb{N}
$ such that $S_{b}(x_{n},x_{n},x_{m})<\varepsilon $ for each $n,m\geq n_{0}$%
. }

\item {\small The $S_{b}$-metric space $(X,S_{b})$ is said to be complete if
every Cauchy sequence is convergent. }
\end{enumerate}

{\small Now we investigate some relationships between $S_{b}$-metric and
some other metrics. The relationship between a metric and an $S$-metric are
given in \cite{Hieu-2015} as follows: }

{\small \textbf{Lemma 2.1.} \label{lem3} \cite{Hieu-2015} Let $(X,d)$ be a
metric space. Then the following properties are satisfied: }

\begin{enumerate}
\item {\small $S_{d}(x,y,z)=d(x,z)+d(y,z)$ for all $x,y,z\in X$ is an $S$%
-metric on $X$. }

\item {\small $x_{n}\rightarrow x$ in $(X,d)$ if and only if $%
x_{n}\rightarrow x$ in $(X,S_{d})$. }

\item {\small $\{x_{n}\}$ is Cauchy in $(X,d)$ if and only if $\{x_{n}\}$ is
Cauchy in $(X,S_{d})$. }

\item {\small $(X,d)$ is complete if and only if $(X,S_{d})$ is complete. }
\end{enumerate}

{\small Since every $S$-metric is an $S_{b}$-metric, using Lemma $2.1$ we
say that an $S_{b}$-metric generated by a metric $d$ is defined as follows:%
\begin{equation*}
S_{b}^{d}(x,y,z)=b[d(x,z)+d(y,z)]\text{,}
\end{equation*}%
for all $x,y,z\in X$ with $b\geq 1$. But there exists an $S_{b}$-metric
which is not generated by any metric as seen in the following example. }

{\small \textbf{Example 2.5.} \label{exm3} Let $X=%
\mathbb{R}
$. We consider the function $S:X\times X\times X\rightarrow \lbrack 0,\infty
)$ given in \cite{Nihal-2015-2} as follows:%
\begin{equation*}
S(x,y,z)=\left\vert x-z\right\vert +\left\vert x+z-2y\right\vert \text{,}
\end{equation*}%
for all $x,y,z\in
\mathbb{R}
$. Then $(X,S)$ is an $S$-metric space. Hence $(X,S)$ is an $S_{b}$-metric
space with $b=1$. This metric is not generated by any metric $d$. }

{\small In the following lemmas, we show that the relationships between $b$%
-metric and $S_{b}$-metric. }

{\small \textbf{Lemma 2.2.} \label{lem4} Let $(X,S_{b})$ be an $S_{b}$%
-metric space, $S_{b}$ be a symmetric $S_{b}$-metric with $b\geq 1$ and the
function $d:X\times X\rightarrow \lbrack 0,\infty )$ be defined by%
\begin{equation*}
d(x,y)=S_{b}(x,x,y)\text{,}
\end{equation*}%
for all $x,y\in X$. Then $d$ is a $b$-metric on $X$. }

{\small
\begin{proof}
It can be easily seen that the conditions $(b1)$ and $(b2)$ are satisfied.
Now we show that the condition $(b3)$ is satisfied. Using the inequality $%
(S_{b}2)$ we have%
\begin{eqnarray*}
d(x,y) &=&S_{b}(x,x,y)\leq b[2S_{b}(x,x,z)+S_{b}(y,y,z)] \\
&=&2bS_{b}(x,x,z)+bS_{b}(y,y,z)
\end{eqnarray*}%
and%
\begin{eqnarray*}
d(x,y) &=&S_{b}(y,y,x)\leq b[2S_{b}(y,y,z)+S_{b}(x,x,z)] \\
&=&2bS_{b}(y,y,z)+bS_{b}(x,x,z)\text{.}
\end{eqnarray*}%
Hence we obtain%
\begin{equation*}
d(x,y)\leq \dfrac{3b}{2}[d(x,z)+d(y,z)]\text{,}
\end{equation*}%
for all $x,y\in X$. Then $d$ is a $b$-metric on $X$ with $\dfrac{3b}{2}$.
\end{proof}
}

{\small \textbf{Lemma 2.3.} \label{lem5} Let $(X,d)$ be a $b$-metric space
with $b\geq 1$ and the function $S_{b}:X\times X\times X\rightarrow \lbrack
0,\infty )$ be defined by%
\begin{equation*}
S_{b}(x,y,z)=d(x,z)+d(y,z)\text{,}
\end{equation*}%
for all $x,y,z\in X$. Then $S_{b}$ is an $S_{b}$-metric on $X$. }

{\small
\begin{proof}
It can be easily verify that the conditions $(S_{b}1)$ is satisfied. We
prove that the condition $(S_{b}2)$ is satisfied. Using the inequality $(b3)$
we get%
\begin{eqnarray*}
S_{b}(x,y,z) &=&d(x,z)+d(y,z) \\
&\leq &b[d(x,a)+d(a,z)]+b[d(y,a)+d(a,z)] \\
&=&bd(x,a)+2bd(a,z)+bd(y,a) \\
&\leq &2bd(x,a)+2bd(y,a)+2bd(a,z) \\
&=&b[S_{b}(x,x,a)+S_{b}(y,y,a)+S_{b}(z,z,a)]\text{,}
\end{eqnarray*}%
for all $x,y,z\in X$. Then $S_{b}$ is an $S_{b}$-metric on $X$ with $b$.
\end{proof}}

{\small Now we give the following example to show that there exists an $%
S_{b} $-metric which is not generated by any $b$-metric. }

{\small \textbf{Example 2.6.} \label{exm9} Let $X=%
\mathbb{R}
$ and define the function $S_{b}:X\times X\times X\rightarrow \lbrack
0,\infty )$%
\begin{equation*}
S_{b}(x,y,z)=b\left( \left\vert x-z\right\vert +\left\vert x+z-2y\right\vert
\right) \text{,}
\end{equation*}%
for all $x,y,z\in
\mathbb{R}
$, where $b\geq 1$. Then $(%
\mathbb{R}
,S_{b})$ is an $S_{b}$-metric space. Now we show that there does not exist
any $b$-metric $d$ which generates this $S_{b}$-metric. Conversely, assume
that there exists a $b$-metric $d$ such that%
\begin{equation*}
S_{b}(x,y,z)=d(x,z)+d(y,z)\text{,}
\end{equation*}%
for all $x,y,z\in
\mathbb{R}
$. Then we get%
\begin{equation*}
S_{b}(x,x,z)=2d(x,z)=2b\left\vert x-z\right\vert \text{ and }%
d(x,z)=b\left\vert x-z\right\vert
\end{equation*}%
and%
\begin{equation*}
S_{b}(y,y,z)=2d(y,z)=2b\left\vert y-z\right\vert \text{ and }%
d(y,z)=b\left\vert y-z\right\vert \text{,}
\end{equation*}%
for all $x,y,z\in
\mathbb{R}
$. Therefore we obtain%
\begin{equation*}
b\left( \left\vert x-z\right\vert +\left\vert x+z-2y\right\vert \right)
=b\left\vert x-z\right\vert +b\left\vert y-z\right\vert \text{,}
\end{equation*}%
which is a contradiction. Consequently, the $S_{b}$-metric can not be
generated by any $b$-metric. }

{\small \textbf{Remark 2.1.} \label{rem2} Notice that the class of all $S$%
-metrics and the class of all $G$-metrics are distinct \cite{Dung-2014}.
Since every $S$-metric is an $S_{b}$-metric and evey $G$-metric is a $G_{b}$%
-metric then the class of all $S_{b}$-metrics and the class of all $G_{b}$%
-metrics are distinct. \vspace{0.2in}\newline
\setcounter{equation}{0} \renewcommand{\theequation}{3.\arabic{equation}}
\noindent \textbf{4. Some fixed point results} \vspace{0.2in}\newline
\indent\label{sec:3} In this section we prove the Banach's contraction
principle on complete $S_{b}$-metric spaces. Then we give a generalization
of this principle. We use the following lemma.}

{\small \textbf{Lemma 3.1.} \label{lem7} \cite{Sedghi-2016} Let $(X,S_{b})$
be an $S_{b}$-metric space with $b\geq 1$, then we have
\begin{equation*}
S_{b}(x,x,y)\leq bS_{b}(y,y,x)\text{ and }S_{b}(y,y,x)\leq bS_{b}(x,x,y).%
\text{ }
\end{equation*}%
}

{\small \textbf{Theorem 3.1.} \label{thm1} Let $(X,S_{b})$ be a complete $%
S_{b}$-metric space with $b\geq 1$ and $T:X\rightarrow X$ be a self-mapping
satisfying%
\begin{equation}
S_{b}(Tx,Tx,Ty)\leq hS_{b}(x,x,y)\text{,}  \label{eqn1}
\end{equation}%
for all $x,y,z\in X$, where $0\leq h<\dfrac{1}{b^{2}}$. Then $T$ has a fixed
point $x$ in $X$. }

{\small
\begin{proof}
Let $T$ satisfies the inequality $($\ref{eqn1}$)$ and $x_{0}\in X$. Then we
define the sequence $\{x_{n}\}$ by $x_{n}=T^{n}x_{0}$. Using the inequality $%
($\ref{eqn1}$)$ and mathematical induction we obtain%
\begin{equation}
S_{b}(x_{n},x_{n},x_{n+1})\leq h^{n}S_{b}(x_{0},x_{0},x_{1})\text{.}
\label{eqn2}
\end{equation}%
Since the conditions $(S_{b}2)$ and $($\ref{eqn2}$)$ are satisfied for all $%
n,m\in
\mathbb{N}
$ with $m>n$, using Lemma 3.1 we get%
\begin{eqnarray*}
S_{b}(x_{n},x_{n},x_{m}) &\leq
&b[2S_{b}(x_{n},x_{n},x_{n+1})+S_{b}(x_{m},x_{m},x_{n+1})] \\
...
\\
&\leq &\dfrac{2bh^{n}}{1-b^{2}h}S_{b}(x_{0},x_{0},x_{1})\text{.}
\end{eqnarray*}%
Since $h\in \lbrack 0,\dfrac{1}{b^{2}})$, where $b\geq 1$, taking limit for $%
n\rightarrow \infty $ then we obtain $S_{b}(x_{n},x_{n},x_{m})\rightarrow 0$
and so $\{x_{n}\}$ is a Cauchy sequence. Since $X$ is complete $S_{b}$%
-metric space there exists $x\in X$ with $\underset{n\rightarrow \infty }{%
\lim }x_{n}=x$.

Assume that $Tx\neq x$. Using the inequality $($\ref{eqn1}$)$ we have%
\begin{equation*}
S_{b}(Tx,Tx,x_{n+1})\leq hS_{b}(x,x,x_{n})\text{.}
\end{equation*}%
If we take limit for $n\rightarrow \infty $ we get a contradiction as
follows:%
\begin{equation*}
S_{b}(Tx,Tx,x)\leq hS_{b}(x,x,x)\text{.}
\end{equation*}%
Hence $Tx=x$. Now we show that the fixed point $x$ is unique. Suppose that $%
Tx=x$, $Ty=y$ and $x\neq y$. Using the inequality $($\ref{eqn1}$)$ we have%
\begin{equation*}
S_{b}(Tx,Tx,Ty)=S_{b}(x,x,y)\leq hS_{b}(x,x,y)\text{.}
\end{equation*}%
We obtain $x=y$ since $h\in \lbrack 0,\dfrac{1}{b^{2}})$. Consequently $x$ is a
unique fixed point of the self-mapping $T$.
\end{proof}}

{\small \textbf{Corollary 3.1.} \label{cor1} Let $(X,S_{b})$ be a complete $%
S_{b}$-metric space with $b\geq 1$, $S_{b}$ be symmetric and $T:X\rightarrow
X$ be a self-mapping satisfying the inequality (\ref{eqn1}) for all $%
x,y,z\in X$, where $0\leq h<\dfrac{1}{b}$. Then $T$ has a fixed point $x$ in
$X$. }

{\small \textbf{Example 3.1.} \label{exm4} Let $X=%
\mathbb{R}
$ and consider the $S_{b}$-metric defined in Example $2.1$ as follows:%
\begin{equation*}
S_{b}(x,y,z)=\dfrac{1}{16}(\left\vert x-y\right\vert +\left\vert
y-z\right\vert +\left\vert x-z\right\vert )^{2}\text{,}
\end{equation*}%
for all $x,y,z\in
\mathbb{R}
$ with $b=4$. If we define the self-mapping $T$ of $%
\mathbb{R}
$ as%
\begin{equation*}
Tx=\dfrac{x}{6}\text{,}
\end{equation*}%
for all $x\in
\mathbb{R}
$ then $T$ satisfies the condition of the Banach's contraction principle.
Indeed we get%
\begin{equation*}
S_{b}(Tx,Tx,Ty)=\dfrac{\left\vert x-y\right\vert ^{2}}{144}\leq
hS_{b}(x,x,y)=\dfrac{\left\vert x-y\right\vert ^{2}}{72}\text{,}
\end{equation*}%
for all $x\in
\mathbb{R}
$ and $h=\dfrac{1}{18}$. Hence $T$ has a unique fixed point $x=0$ in $%
\mathbb{R}
$. }

{\small Now we give the following theorem as a generalization of the
Banach's contraction principle on complete $S_{b}$-metric spaces. }

{\small \textbf{Theorem 3.2.} \label{thm2} Let $(X,S_{b})$ be a complete $%
S_{b}$-metric space with $b\geq 1$ and $T$ be a self-mapping of $X$
satisfying the following condition: }

{\small There exist real numbers $\alpha _{1}$, $\alpha _{2}$ satisfying $%
\alpha _{1}+\left( 2b^{2}+b\right) \alpha _{2}<1$ with $\alpha _{1}$, $%
\alpha _{2}\geq 0$ such that%
\begin{eqnarray}
S_{b}(Tx,Tx,Ty) &\leq &\alpha _{1}S_{b}(x,x,y)+\alpha _{2}\max
\{S_{b}(Tx,Tx,x),  \label{eqn3} \\
&&S_{b}(Tx,Tx,y),S_{b}(Ty,Ty,y),S_{b}(Ty,Ty,x)\}\text{,}  \notag
\end{eqnarray}%
for all $x,y\in X$. Then $T$ has a unique fixed point $x$ in $X$. }

{\small
\begin{proof}
Let $x_{0}\in X$ and the sequence $\{x_{n}\}$ be defined as follows:%
\begin{equation*}
Tx_{0}=x_{1}\text{, }Tx_{1}=x_{2}\text{,}\cdots \text{, }Tx_{n}=x_{n+1}\text{%
, }\cdots \text{.}
\end{equation*}%
Assume that $x_{n}\neq x_{n+1}$ for all $n$. Using the condition $($\ref%
{eqn3}$)$ we get%
\begin{eqnarray}
S_{b}(x_{n},x_{n},x_{n+1}) &=&S_{b}(Tx_{n-1},Tx_{n-1},Tx_{n})\leq \alpha
_{1}S_{b}(x_{n-1},x_{n-1},x_{n})  \label{eqn4} \\
&&+\alpha _{2}\max \{S_{b}(x_{n},x_{n},x_{n-1}),S_{b}(x_{n},x_{n},x_{n}),
\notag \\
&&S_{b}(x_{n+1},x_{n+1},x_{n}),S_{b}(x_{n+1},x_{n+1},x_{n-1})\}  \notag \\
&=&\alpha _{1}S_{b}(x_{n-1},x_{n-1},x_{n})+\alpha _{2}\max
\{S_{b}(x_{n},x_{n},x_{n-1}),  \notag \\
&&S_{b}(x_{n+1},x_{n+1},x_{n}),S_{b}(x_{n+1},x_{n+1},x_{n-1})\}\text{.}
\notag
\end{eqnarray}%
By the condition $(S_{b}2)$ we have%
\begin{equation}
S_{b}(x_{n+1},x_{n+1},x_{n-1})\leq
b[2S_{b}(x_{n+1},x_{n+1},x_{n})+S_{b}(x_{n-1},x_{n-1},x_{n})]\text{.}
\label{eqn5}
\end{equation}%
Using the conditions $($\ref{eqn4}$)$, $($\ref{eqn5}$)$ and Lemma 3.1 we obtain%
\begin{eqnarray*}
S_{b}(x_{n},x_{n},x_{n+1}) &\leq &\alpha
_{1}S_{b}(x_{n-1},x_{n-1},x_{n})+\alpha _{2}\max
\{S_{b}(x_{n},x_{n},x_{n-1}), \\
&&S_{b}(x_{n+1},x_{n+1},x_{n}),2bS_{b}(x_{n+1},x_{n+1},x_{n})+bS_{b}(x_{n-1},x_{n-1},x_{n})
\\
&\leq &\alpha _{1}S_{b}(x_{n-1},x_{n-1},x_{n})+2b\alpha
_{2}S_{b}(x_{n+1},x_{n+1},x_{n}) \\
&&+b\alpha _{2}S_{b}(x_{n-1},x_{n-1},x_{n})
\end{eqnarray*}%
and so%
\begin{equation*}
(1-2b^{2}\alpha _{2})S_{b}(x_{n},x_{n},x_{n+1})\leq (\alpha _{1}+b\alpha
_{2})S_{b}(x_{n-1},x_{n-1},x_{n})\text{,}
\end{equation*}%
which implies%
\begin{equation}
S_{b}(x_{n},x_{n},x_{n+1})\leq \dfrac{\alpha _{1}+b\alpha _{2}}{1-2b^{2}\alpha
_{2}}S_{b}(x_{n-1},x_{n-1},x_{n})\text{.}  \label{eqn6}
\end{equation}%
Let $d=\dfrac{\alpha _{1}+b\alpha _{2}}{1-2b^{2}\alpha _{2}}$. Then $d<1$ since $%
\alpha _{1}+\left( 2b^{2}+b\right) \alpha _{2}<1$. Notice that $1-2b^{2}\alpha _{2}\neq 0$ since $%
0\leq \alpha _{2}<\dfrac{1}{2b^{2}+b}$. Now repeating this process in the
inequality $($\ref{eqn6}$)$ we get%
\begin{equation}
S_{b}(x_{n},x_{n},x_{n+1})\leq d^{n}S_{b}(x_{0},x_{0},x_{1})\text{.}
\label{eqn7}
\end{equation}%
We show that the sequence $\{x_{n}\}$ is Cauchy. Then for all $n,m\in
\mathbb{N}
$, $m>n$, using the conditions $($\ref{eqn7}$)$ and $(S_{b}2)$ we obtain%
\begin{eqnarray*}
S_{b}(x_{n},x_{n},x_{m}) \leq &\dfrac{2bd^{n}}{1-b^{2}d}S_{b}(x_{0},x_{0},x_{1})\text{.}  \notag
\end{eqnarray*}%
We have $\underset{n,m\rightarrow \infty }{\lim }S_{b}(x_{n},x_{n},x_{m})=0$
by the above inequality and so $\{x_{n}\}$ is a Cauchy sequence. By
the completeness hypothesis, there exists $x\in X$ such that $\{x_{n}\}$
converges to $x$. Suppose that $Tx\neq x$. Then we have%
\begin{eqnarray*}
S_{b}(x_{n},x_{n},Tx) &=&S_{b}(Tx_{n-1},Tx_{n-1},Tx) \\
&\leq &\alpha _{1}S_{b}(x_{n-1},x_{n-1},x)+\alpha _{2}\max
\{S_{b}(x_{n},x_{n},x_{n-1}), \\
&&S_{b}(x_{n},x_{n},x),S_{b}(Tx,Tx,x),S_{b}(Tx,Tx,x_{n-1})\}
\end{eqnarray*}%
and so taking limit for $n\rightarrow \infty $ and using Lemma 3.1 we get%
\begin{equation*}
S_{b}(x,x,Tx)\leq \alpha _{2}S_{b}(Tx,Tx,x)\leq \alpha _{2}bS_{b}(x,x,Tx)\text{,}
\end{equation*}%
which implies $S_{b}(Tx,Tx,x)=0$ and $Tx=x$ since $0\leq \alpha _{2}<\dfrac{1%
}{2b^{2}+b}$.

Finally we show that the fixed point $x$ is unique. To do this we assume
that $x\neq y$ such that $Tx=x$ and $Ty=y$. Using the inequality $(\ref{eqn3}%
)$ and Lemma 3.1 we have%
\begin{eqnarray*}
S_{b}(Tx,Tx,Ty) &=&S_{b}(x,x,y)\leq \alpha _{1}S_{b}(x,x,y)+\alpha _{2}\max
\{S_{b}(x,x,x), \\
&&S_{b}(x,x,y),S_{b}(y,y,y),S_{b}(y,y,x)\}\text{,}
\end{eqnarray*}%
which implies $x=y$ since $\alpha _{1}+b\alpha _{2}<1$. Then the proof is
completed.
\end{proof}
}

{\small \textbf{Corollary 3.2.} \label{cor2} Let $(X,S_{b})$ be a complete $%
S_{b}$-metric space with $b\geq 1$, $S_{b}$ be symmetric and $T$ be a
self-mapping of $X$ satisfying the following condition: }

{\small There exist real numbers $\alpha _{1}$, $\alpha _{2}$ satisfying $%
\alpha _{1}+3b\alpha _{2}<1$ with $\alpha _{1}$, $\alpha _{2}\geq 0$ such
that%
\begin{eqnarray*}
S_{b}(Tx,Tx,Ty) &\leq &\alpha _{1}S_{b}(x,x,y)+\alpha _{2}\max
\{S_{b}(Tx,Tx,x), \\
&&S_{b}(Tx,Tx,y),S_{b}(Ty,Ty,y),S_{b}(Ty,Ty,x)\}\text{,}  \notag
\end{eqnarray*}%
for all $x,y\in X$. Then $T$ has a unique fixed point $x$ in $X$. }

{\small \textbf{Remark 3.1.} \label{rem3} If we take $b=1$ in Theorem $3.2$
then we obtain Theorem $1 $ in \cite{Nihal-2015-3}. }

{\small \textbf{Remark 3.2.} \label{rem4} We note that Theorem $3.2$ is a
generalization of the Banach's contraction principle on $S_{b}$-metric
spaces. Indeed, if we take $\alpha _{1}<\dfrac{1}{b^{2}}$ and $\alpha _{2}=0$
in Theorem $3.2$ we obtain the Banach's contraction principle. }

{\small Now we give an example of a self-mapping satisfying the conditions
of Theorem $3.2$ such that the condition of the Banach's contraction
principle is not satisfied. }

{\small \textbf{Example 3.2.} \label{exm5} We consider the $S$-metric space $%
(%
\mathbb{R}
,S)$ with%
\begin{equation*}
S(x,y,z)=\left\vert x-z\right\vert +\left\vert x+z-2y\right\vert \text{,}
\end{equation*}%
for all $x,y,z\in
\mathbb{R}
$ given in \cite{Nihal-2015-2} and the self-mapping $T$ of $%
\mathbb{R}
$ as%
\begin{equation*}
Tx=\left\{
\begin{array}{ll}
x+50 & \text{if }\left\vert x-1\right\vert =1 \\
45 & \text{if }\left\vert x-1\right\vert \neq 1%
\end{array}%
\right. \text{,}
\end{equation*}%
for all $x\in
\mathbb{R}
$ defined in \cite{Nihal-2015-3}. Since every $S$-metric space is an $S_{b}$%
-metric space, $(%
\mathbb{R}
,S)$ is an $S_{b}$-metric space with $b=1$. Then the inequality $($\ref{eqn3}%
$)$ is satisfied for $\alpha _{1}=0$ and $\alpha _{2}=\dfrac{1}{5}$. Then $T$
has a unique fixed point $x=45$ by Theorem $3.2$. But $T$ does not satisfy
the condition of the Banach's contraction principle since for $x=1$, $y=0$
we get%
\begin{equation*}
S(Tx,Tx,Ty)=10\leq hS(x,x,y)=2h\text{,}
\end{equation*}%
which is a contradiction with $h<1$. \vspace{0.2in}\newline
\setcounter{equation}{0} \renewcommand{\theequation}{4.\arabic{equation}}
\noindent \textbf{5. An application of the Banach's contraction to linear
equations} \vspace{0.2in}\newline
\indent\label{sec:4} In this section we give an application of the Banach's
contraction principle on $S_{b}$-metric spaces to linear equations. To do
this we consider the $S_{b}$-metric space generated by%
\begin{equation*}
d_{1}(x,y)=\sum\limits_{i=1}^{n}\left\vert x_{i}-y_{i}\right\vert \text{,}
\end{equation*}%
for all $x,y\in
\mathbb{R}
^{n}$. We note that the symmetry condition (\ref{symmetry}) is not necessary
in the following example. }

{\small \textbf{Example 4.1.} \label{exm6} Let $X=%
\mathbb{R}
^{n}$ be an $S_{b}$-metric space with the $S_{b}$-metric defined by%
\begin{equation*}
S_{1}(x,y,z)=\sum\limits_{i=1}^{n}\left\vert x_{i}-z_{i}\right\vert
+\sum\limits_{i=1}^{n}\left\vert y_{i}-z_{i}\right\vert \text{,}
\end{equation*}%
for all $x,y,z\in
\mathbb{R}
^{n}$, where $b=1$. If%
\begin{equation*}
\underset{1\leq j\leq n}{\max }\sum\limits_{i=1}^{n}\left\vert
a_{ij}\right\vert \leq h<1\text{,}
\end{equation*}%
then the linear system%
\begin{eqnarray}
a_{11}x_{1}+a_{12}x_{2}+\cdots +a_{1n}x_{n} &=&b_{1}  \label{eqn9} \\
a_{21}x_{1}+a_{22}x_{2}+\cdots +a_{2n}x_{n} &=&b_{2}  \notag \\
&&\vdots  \notag \\
a_{n1}x_{1}+a_{n2}x_{2}+\cdots +a_{nn}x_{n} &=&b_{n}  \notag
\end{eqnarray}%
has a unique solution. Let the self-mapping $T:X\rightarrow X$ be defined by%
\begin{equation*}
Tx=Ax+b\text{,}
\end{equation*}%
where $x,b\in
\mathbb{R}
^{n}$ and%
\begin{equation*}
\left(
\begin{array}{cccc}
a_{11} & a_{12} & \cdots & a_{1n} \\
a_{21} & a_{22} & \cdots & a_{2n} \\
\vdots & \vdots & \ldots & \vdots \\
a_{n1} & a_{n2} & \cdots & a_{nn}%
\end{array}%
\right) \text{.}
\end{equation*}%
Now we show that the self-mapping satisfies the contraction of the Banach's
contraction principle. For $x,y\in
\mathbb{R}
^{n}$ we get%
\begin{eqnarray*}
S_{1}(Tx,Tx,Ty) &=&2\sum\limits_{i=1}^{n}\left\vert
\sum\limits_{j=1}^{n}a_{ij}(x_{j}-y_{j})\right\vert \leq
2\sum\limits_{i=1}^{n}\sum\limits_{j=1}^{n}\left\vert a_{ij}\right\vert
\left\vert x_{j}-y_{j}\right\vert \\
&=&2\sum\limits_{j=1}^{n}\sum\limits_{i=1}^{n}\left\vert a_{ij}\right\vert
\left\vert x_{j}-y_{j}\right\vert =\sum\limits_{j=1}^{n}2\left\vert
x_{j}-y_{j}\right\vert \sum\limits_{i=1}^{n}\left\vert a_{ij}\right\vert \\
&\leq &\left( \underset{1\leq j\leq n}{\max }\sum\limits_{i=1}^{n}\left\vert
a_{ij}\right\vert \right) \sum\limits_{j=1}^{n}2\left\vert
x_{j}-y_{j}\right\vert \leq hS_{1}(x,x,y)\text{.}
\end{eqnarray*}%
Then $T$ satisfy the contraction of the Banach's contraction principle.
Using Theorem $3.2$ the linear equations system $($\ref{eqn9}$)$ has a
unique solution. \vspace{0.2in}\newline
{\footnotesize {\ }}}

{\small {\footnotesize \ }}

\end{document}